\documentclass[12pt]{amsart}

\title{}
\author{}

\usepackage{fullpage, amsfonts, amsmath, amscd, amsthm,stmaryrd}
\usepackage{amssymb}
\usepackage{graphics}
\usepackage{graphicx}
\usepackage{amscd}

\newtheorem{theorem}{Theorem}[section]
\newtheorem{lemma}[theorem]{Lemma}

\newtheorem*{thm*}{Theorem}

\theoremstyle{definition}

\newcommand{\A}{\mathbb{A}}
\newcommand{\Z}{\mathbb{Z}}
\newcommand{\R}{\mathbb{R}}

\newcommand\goth[1]{\mathfrak{#1}}
\newcommand{\s}{\goth{s}}

\title{Perturbed Floer Homology of some fibered three manifolds II}
\author{Zhongtao Wu}

\begin{document}

\maketitle

\begin{abstract}
Modifying the method of \cite{ztwu}, we compute the perturbed $HF^+$ for some special classes of fibered three manifolds in the second highest spin$^c$-structures $S_{g-2}$.  The special classes considered in this paper include the mapping tori of Dehn twists along a single non-separating curve and along a transverse pair of curves.

\end{abstract}

\section{Introduction}

Following \cite{ztwu}, where the perturbed Heegaard Floer homology was defined and computed for the product three manifolds $\Sigma_g \times S^1$, we aim to modify our method for the computations of certain general fibered three manifolds.  More precisely, we treat each fibered three manifold as a mapping torus $M(\phi)$ for some orientation-preserving diffeomorphism $\phi: \Sigma_g \longrightarrow \Sigma_g$, and decompose $\phi$ into products of {\it Dehn twists}.  The cases studied here consist of those Dehn twists along a single non-separating curve, and those along a transverse pair of curves.  

Fibered three manifolds admit certain particularly simple ``special Heegaard Diagrams'', first introduced by Ozsv\'ath and Szab\'o in \cite {OSzCont}, where a genus $2g+1$ Heegaard Diagrams was constructed for each mapping torus $M(\phi)$.  Let $S_k \subset Spin^c(M(\phi))$ denote the set of spin$^c$-structures $\s$ with $\langle c_1(\s),[\Sigma_g]\rangle =2k.$  We will focus on the computation of the homology in the set of spin$^c$-structures $S_{g-2}$, for reasons to explain shortly. 
The details for the backgrounds of the ``special Heegaard Diagram'' are reviewed in section 2.  

Section 3 through section 7 are dealing with various special classes of fibered three manifolds.  All these cases are approached in a similar manner: We write down a special Heegaard diamgram for each manifold, and find all the generators in $S_{g-2}$.  The Euler characteristic is subsequently computed in each case, being made possible due to Lemma \ref{l1} and \ref{l2}, where it is identified with the {\it Lefschetz number}  $L(\phi):= \sum_i (-1)^i trace(\phi_*: H_i(M)\rightarrow H_i(M))$, and hence reduces the computation to a simple matter of linear algebra.  The number of the generators and the Euler characteristic are found to be equal in each case, so we follow with the argument of \cite{ztwu}.  In the end, the homology is found to be equal to the Euler characteristic in each spin$^c$-structure.  

Results of a similar nature have been obtained by various other people.  Seidel \cite{Se} considered the {\it symplectic Floer homology} of surface symplectomorphisms, calculating it for arbitrary compositions of Dehn twists along a disjoint collection curves.  Eftekhary \cite{E} generalized Seidel's work to Dehn twists along two disjoint forests.  It is then Cotton-Clay who achieved a vast generalization to include all pseudo-Anosov mapping class and reducible mapping class.  Jabuka and Mark \cite{JM}, on the other hand, computed the unperturbed results for the {\it Heegaard Floer homology}.  Their results agreed with the previously mentioned papers wherever applicable, presenting a strong piece of favorable evidence for the conjectural existence of the isomorphisms between all versions of Floer homologies.  Very recently, Taubes \cite{T} claimed a proof for the equivalence between {\it Seiberg-Witten Floer cohomology} and the {\it embedded contact homology}, though the equivalences of other versions of Floer homology are not established yet.  Our paper is largely motivated by this, for the spin$^c$-structures in $S_{g-2}$ are the relevant parts of the Heegaard Floer homology in the conjecture.

For the effects of the perturbations in Heegaard Floer theory are still poorly understood at the time being, the results in our paper are also intended as interesting examples for studying and understanding the relation.  We find our perturbed homologies in section 4 agree with Jabuka and Mark's unperturbed results \cite{JM}, while they are strictly smaller than the latter counterparts in section 3 and 5.  It was also mentioned that all the homologies in this paper are identical to the corresponding Euler characteristics.  Though it is perhaps too bold to conjecture such a phenomenon occurs for all (fibered) manifolds, we strongly believe it should hold for a much larger class of manifolds than those being discovered in this article.  It may worth the effort to explore further in this direction.

\subsection*{Acknowledgment.}  I am much obliged to my advisor, Zolt\'an  Szab\'o, for his continual encouragements and preparing me with the necessary backgrounds.  I am also grateful to Joshua Greene and Yi Ni for helpful discussions at various points.

\section{Review of the special Heegaard Diagram}

In this section, we review the special Heegaard Diagram construction by Ozsv\'ath and Szab\'o \cite[section 3]{OSzCont}.  Figure \ref{Sigma2} is the Heegaard Diagram for $\Sigma_g\times S^1$ used in \cite{ztwu}, consisting of two 2g-gons with standard identifications of edges and two punctured holes.  They represent two genus $g$ surfaces, joined together through the pairs of holes to make a genus $2g+1$ surface.  All the $\alpha$'s and $\beta$'s curves are drawed along with their intersection points marked.  We list all the interesting properties of the Heegaard Diagram:

\begin{figure}
\begin{center}
\includegraphics[width=4in]{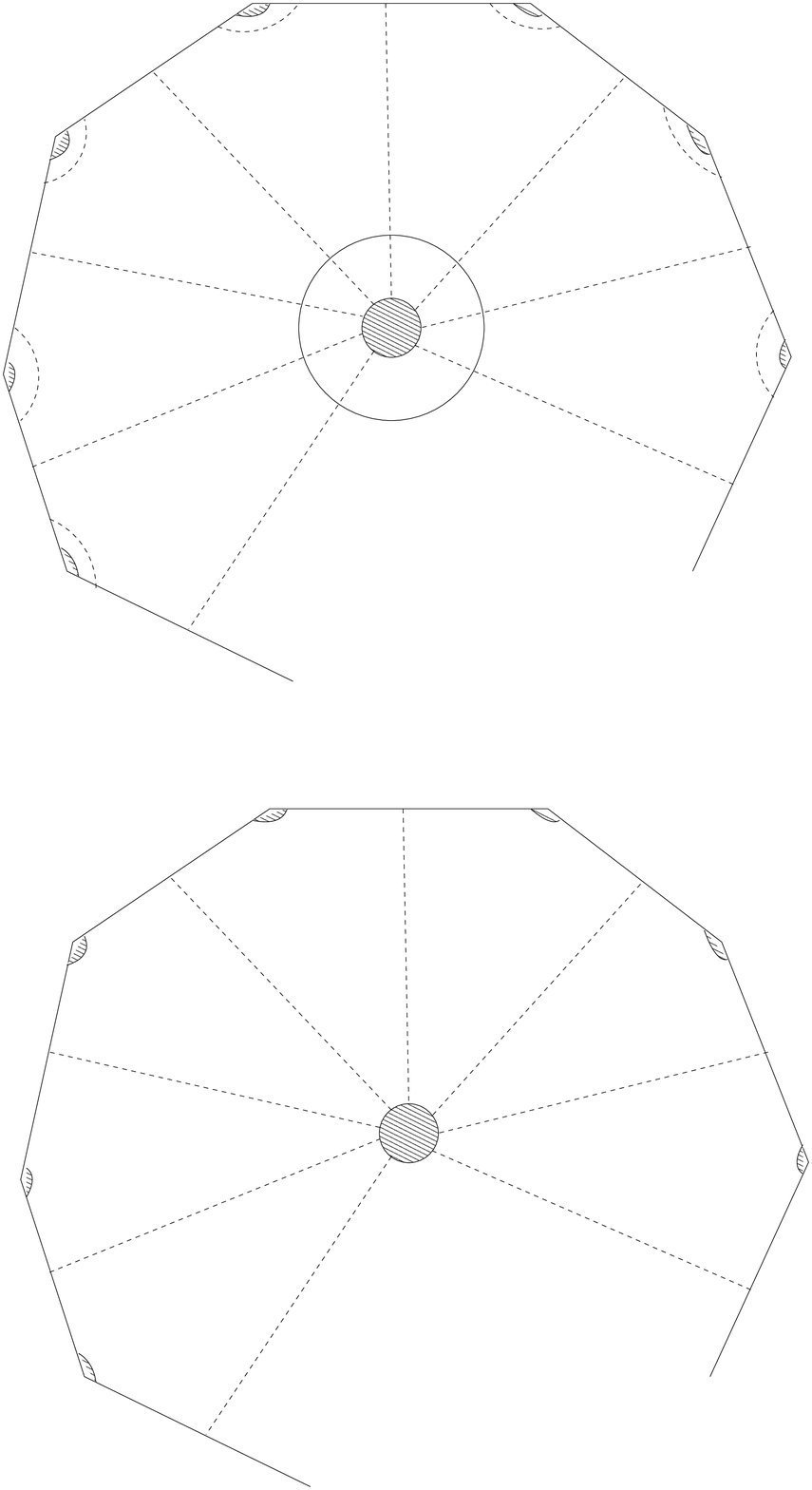}
\setlength{\unitlength}{0.0112in}
\put(-15,100){$\alpha_2$}
\put(-27,210){$\alpha_1$}
\put(-80,285){$\alpha_2$}
\put(-210,310){$\alpha_1$}
\put(-190,305){$R_1$}
\put(-75,265){$R_2$}
\put(-20,190){$R_1$}
\put(-20,80){$R_2$}
\put(-320,270){$R_{2g}$}
\put(-305,285){$\alpha_{2g}$}
\put(-370,210){$\alpha_{2g-1}$}
\put(-350,70){$\alpha_{2g}$}
\put(-310,20){$\alpha_{2g-1}$}
\put(-360,85){$R_{2g}$}
\put(-385,195){$R_{2g-1}$}
\put(-20,460){$\alpha_2$}
\put(-32,570){$\alpha_1$}
\put(-85,645){$\alpha_2$}
\put(-70,625){$L_2$}
\put(-215,670){$\alpha_1$}
\put(-25,550){$L_1$}
\put(-195,665){$L_1$}
\put(-30,440){$L_2$}
\put(-325,630){$L_{2g}$}
\put(-310,645){$\alpha_{2g}$}
\put(-390,530){$\alpha_{2g-1}$}
\put(-375,570){$B_{2g-1}'$}
\put(-390,550){$L_{2g-1}$}
\put(-390,510){$B_{2g-1}$}
\put(-380,470){$B_{2g}$}
\put(-360,420){$B_{2g}'$}
\put(-365,445){$L_{2g}$}
\put(-380,430){$\alpha_{2g}$}
\put(-315,380){$\alpha_{2g-1}$}
\put(-200,600){$\beta_1$}
\put(-140,580){$\beta_2$}
\put(-100,540){$\beta_1$}
\put(-105,480){$\beta_2$}
\put(-255,575){$\beta_{2g}$}
\put(-290,525){$\beta_{2g-1}$}
\put(-240,515){$A_{2g-1}$}
\put(-245,495){$A_{2g}$}
\put(-230,540){$A_{2g-1}'$}
\put(-230,475){$A_{2g}'$}
\put(-205,560){$A_1'$}
\put(-170,550){$A_2$}
\put(-150,525){$A_1$}
\put(-145,495){$A_2'$}
\put(-270,480){$\beta_{2g}$}
\put(-235,440){$\beta_{2g-1}$}
\put(-195,240){$\beta_1$}
\put(-135,220){$\beta_2$}
\put(-95,180){$\beta_1$}
\put(-100,120){$\beta_2$}
\put(-250,215){$\beta_{2g}$}
\put(-265,165){$\beta_{2g-1}$}
\put(-265,120){$\beta_{2g}$}
\put(-230,80){$\beta_{2g-1}$}
\put(-160,480){$\alpha_{2g+1}$}
\put(-50,500){$\beta_{2g+1}$}
\put(-280,500){$z$}
\put(-310,500){$D$}
\put(-310,140){$D'$}
\caption{\label{Sigma2} The special Heegaard Diagram for $\Sigma_g\times S^1$. }
\end{center}

\end{figure}

\begin{itemize}
\item each $\alpha_i \cap \beta_i$ twice, denoted by $L_i$ and $R_i$ respectively, $1\leq i \leq 2g$.  

\item $\alpha_i \cap \beta_j= \emptyset$, when $i\neq j$, $1 \leq i, j \leq 2g$.

\item $\alpha_{2g+1} \cap \beta_i$ twice, denoted by $A_i$ and $A'_i$ respectively, $1\leq i \leq 2g$.

\item $\alpha_i \cap \beta_{2g+1}$ twice, denoted by $B_i$ and $B'_i$ respectively, $1\leq i \leq 2g$. 
              
\end{itemize}

Next, we enumerate all generators in this Heegaard diagram.  We sort them according to their Spin$^c$ structures:
\begin{itemize}
 \item When $k \geq g$, $S_k$ is empty.
\item When $k=g-1$, $S_{g-1}$ consists of a pair of generators: $(A_{2g}, B_{2g}, L_1, L_2,\cdots,L_{2g-1})$ and $(A_{2g-1},B_{2g-1}, L_1,\cdots,L_{2g-2},L_{2g})$. 
\item When $k=g-2$, $S_{g-2}$ consists of $(2g-1)$ pairs of generators: \\
$a_1:=(A_{2g}, B_{2g}, R_1, L_2,\cdots,L_{2g-1})$, \\
$a_2:=(A_{2g}, B_{2g}, L_1, R_2, \cdots,L_{2g-1})$,\\
...\\
$a_{2g-2}=(A_{2g}, B_{2g}, L_1,L_2, \cdots, R_{2g-2})$ \\

and \\
$b_1:=(A_{2g-1},B_{2g-1},R_1, L_2,\cdots, L_{2g})$,\\
$b_2:=(A_{2g-1},B_{2g-1},L_1, R_2,\cdots, L_{2g})$,\\
...\\
$b_{2g-2}=A_{2g-1},B_{2g-1},L_1,L_2 \cdots,R_{2g-2})$\\

and\\
$a_0:=(A_{2g}, B_{2g}, L_1, L_2,\cdots,R_{2g-1})$\\
$b_0:=(A_{2g-1},B_{2g-1},L_1, L_2,\cdots, R_{2g})$.\\

Since $a_0$ and $b_0$ are connected by a disk $D'$ not containing the basepoint $z$, we do not expect them to survive in the homology, and consequently we call them {\it fake generators}.  The remaining $(2g-2)$ pairs, on the other hand, are called {\it essential generators}.
\\ 

\item  When $0<k<g-1$, $S_k$ consists of $\binom{2g-1}{g-1-k}$ pairs of generators: Simply replace $(g-1-k)$ of $L_i$ by $R_i$ in the coordinates of the two generators of $S_{g-1}$.  Among them, $\binom{2g-2}{g-2-k}$ pairs are fake and $\binom{2g-2}{g-1-k}$ pairs are essential.                                                                                                                                          
\end{itemize}

In this way, we organize all the generators of the Heegaard diagram systematically.  Such a schemetic presentation of generators is as well available for a general fibered three manifold $M(\phi)$, and remains relatively simple for $S_{g-2}$. 

Throughout the paper, $g$ is implicitly assumed to be greater than 2.

\section{Multiple Dehn twists along a non-separating curve}

\begin{figure}

\begin{center}
\includegraphics[width=4.5in]{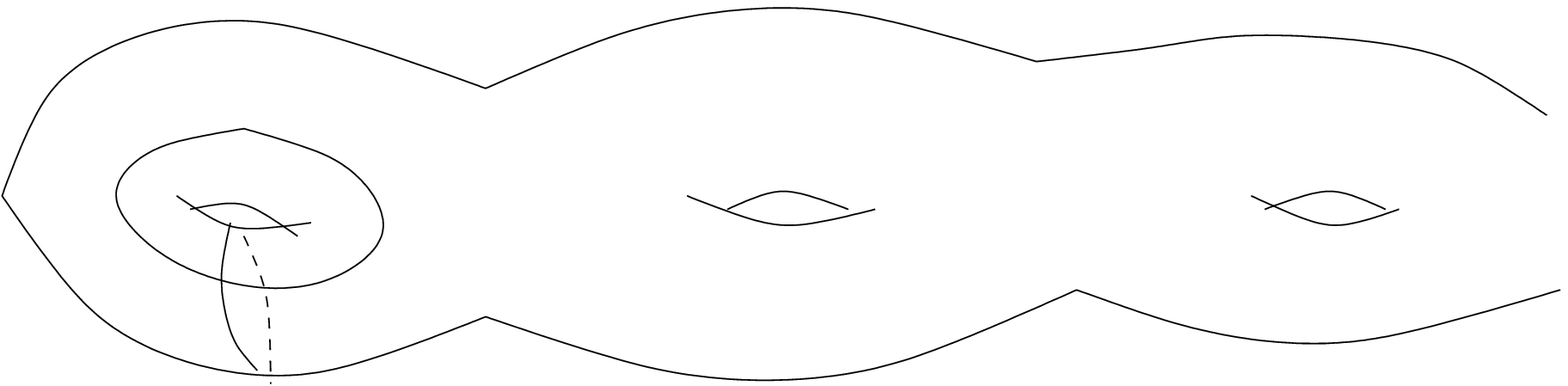}
\setlength{\unitlength}{1.02in}
\put(-3.9,0.1){$\gamma$}
\put(-3.3,0.4){$\delta$}
\caption{\label{gammacurve} } 
\end{center}

\end{figure}

By a standard classification result in surface, any simple non-separating curve can be mapped to the standard position, namely $\gamma$ in Figure \ref{gammacurve}, by a suitable surface automorphism.  Hence, we may assume our manifolds to be $M(t_{\gamma}^n)$ without loss of generality, where $t_\gamma$ denotes the right-handed Dehn twist along $\gamma$.

Next, we draw the special Heegaard diagram of $M(t_{\gamma}^n)$.  In general, for arbitrary $M(\phi)$, the $\alpha$'s and $\beta$'s curves inside the left-hand-side $2g$-gon are the same as that of $\Sigma_g \times S^1$.  As for the right-hand-side $2g$-gon, whereas the $\alpha$'s curves are unaltered, the $\beta$'s curves are twisted according to $\phi$.  Therefore, we would only exhibit the right-hand-side $2g$-gon of the Heegaard diagram, in which all information of the manifolds are encoded.  
 
\begin{figure}

\begin{center}
\includegraphics[width=3in]{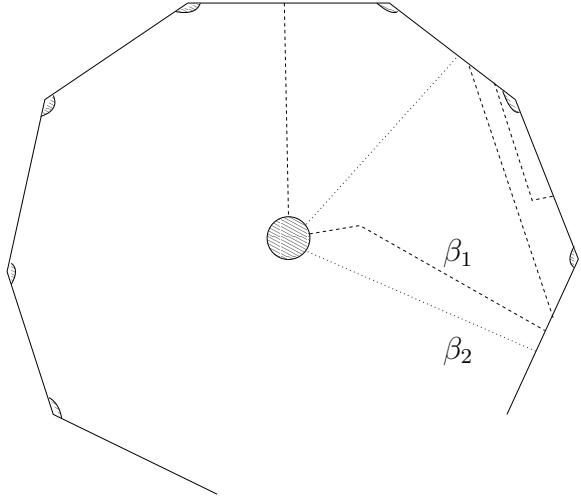}
\setlength{\unitlength}{1.02in}
\put(-0.7,1.2){$\beta_1$}
\put(-0.7,0.7){$\beta_2$}
\caption{\label{HD1} The Heegaard Diagram for $M(t_\gamma^n)$, when $n=2$.  } 
\end{center}

\end{figure}

We proceed to enumerate all the generators in the set of the spin$^c$-structures $S_{g-2}$ in the Heegaard diagram (Figure \ref{HD1}). Observe that a Dehn twist along $\gamma$ does not introduce any new intersections between $\alpha_i$ and $\beta_i$.  Consequently, there is no additional generator, other than the $(2g-1)$ pairs we initially had for $\Sigma_g\times S^1$.

Recall the following identity for the Eular characteristic of $HF^+$ \cite[Proposion 2.3]{ztwu}:

\begin {lemma}\label{l1}
 When $\s$ is a non-torsion Spin$^c$ structure, $HF^+(Y,s;\eta)$ is finitely generated, and the Euler characteristic $$ \chi(HF^+(Y,\s;\eta))=\chi(HF^+(Y,s))=\pm\tau_t(Y,\s),$$ where $\tau_t$ is Turaev's torsion function, with respect to the component $t$ of $H^2(Y;\R)-0$ containing $c_1(\s)$. 

\end {lemma}
 
Turaev's torsion function, derived from some complicated group rings over the CW-complex, is in general rather hard to compute.  In our case for fibered three manifolds though, it is remarkably related to the Lefschetz numbers by the following identity \cite{S},\cite{HL}.

\begin{lemma}\label{l2}

If we denote $\tau_t(M(\phi), k)$ for the sum of all Turaev's torsion functions over the set of the spin$^c$-structures $S_k$, then
 $$\tau_t(M(\phi),k)=L(S^{g-1-k}\phi).$$
where the latter is the Lefschetz number of the induced function of $\phi$ over the symmetric product $S^{g-1-k}\Sigma_g$.  
In particular, for $k=g-2$, $$\tau_t(M(\phi),g-2)=L(\phi).$$ 

\end{lemma}

Hence, applying Lemma \ref{l1} and \ref{l2}, we have $\chi(HF^+ (M(t_{\gamma}^n),\s_{g-2} ; \omega))=L(t_{\gamma}^n)=2-2g$.  Following the argument of \cite[section 4]{ztwu}, we conclude: 
\begin{theorem}

$HF^+ (M(t_{\gamma}^n),S_{g-2} ; \omega)=\A^{2g-2}.$

\end{theorem}
It is interesting to compare our result to those of the unperturbed Heegaard Floer homology \cite{JM} and the symplectic Floer homology \cite {Se}: $HF^+(M(t_{\gamma}),\s_{g-2})\cong \Z^{2g-1}_{(g-1)}\oplus\Z_{(g)}$,  $HF^*(t_{\gamma}^n)\cong H^*(\Sigma_g, \gamma)$.  We find the ranks of our perturbed homology are diminished by 2, a typical phenomenon to anticipate from our previous experiences with computing $\Sigma_g\times S^1$.

We also remark that a generalization is readily attainable for the mapping torus $M(t_{\gamma_1}^{n_1}\cdots t_{\gamma_k}^{n_k})$ of the composition of a sequence of Dehn twists along some mutually disjoint curves $\gamma$'s in ''standard'' positions, thus covering all the cases considered in \cite{Se}.  The techniques and results are exactly the same.

\section{Multiple Dehn twists along a transverse pair of curves, case $m \cdot n <0$}

As before, we may assume the manifolds of our concern to be $M(t_{\gamma}^mt_{\delta}^n)$, without loss of generality, for $\gamma$,$\delta$ in the standard positions of figure \ref{gammacurve}.

\begin{figure}

\begin{center}
\includegraphics[width=3in]{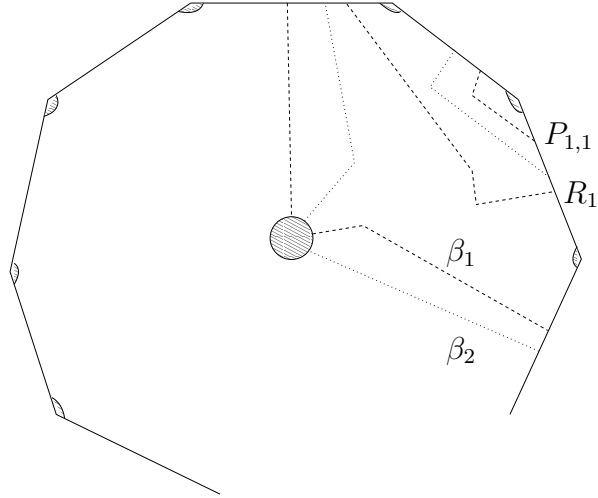}
\setlength{\unitlength}{1.02in}
\put(-0.7,1.2){$\beta_1$}
\put(-0.7,0.7){$\beta_2$}
\put(-0.2,1.8){$P_{1,1}$}
\put(-0.1,1.5){$R_1$}
\caption{\label{HD2} The Heegaard Diagram for $M(t_\gamma^mt_\delta^n)$, when $m=1,n=-1$.  Here, $\beta_1$ is represented by the dashed curve, while $\beta_2$ is represented by the dotted curve.} 
\end{center}

\end{figure}

Consider the case $m \cdot n  <0$ in this section, we have the Heegaard Diagram in Figure \ref{HD2}.  We note extra intersections between $\alpha_1, \alpha_2$ and $\beta_1, \beta_2$, and hence have more generators in $S_{g-2}$.

To be more precise, denote the $|mn|$ extra intersection between $\alpha_1$ and $\beta_1$ by $P_{i,j}$, where $1\leq i\leq |m|$ and $1 \leq j \leq |n|$.  Then, there are $(2g-1+|mn|)$ pairs of generators in $S_{g-2}$, among which $(2g-2+|mn|)$ are essential: \\

$$(A_{2g}, B_{2g}, R_1, L_2,\cdots,L_{2g-1}),$$
$$(A_{2g}, B_{2g}, L_1, R_2, \cdots,L_{2g-1})$$
$$...$$
$$(A_{2g}, B_{2g}, L_1,L_2, \cdots, R_{2g-1})$$
and \\
$$(A_{2g-1},B_{2g-1},R_1, L_2,\cdots, L_{2g})$$
$$(A_{2g-1},B_{2g-1},L_1, R_2,\cdots, L_{2g})$$
$$...$$
$$(A_{2g-1},B_{2g-1},L_1,L_2 \cdots,R_{2g})$$
and \\
$$(A_{2g}, B_{2g}, P_{i,j}, L_2, \cdots, L_{2g-1})$$
$$(A_{2g-1},B_{2g-1}, P_{i,j},L_2, \cdots, L_{2g}).$$

We compute the Lefschetz number of $L(t_{\gamma}^mt_{\delta}^n)$ from its definition
$$ L(\phi)= \sum_i (-1)^i trace(\phi_*: H_i(M)\rightarrow H_i(M))$$

Both $t_\gamma$ and $t_\delta$ act trivially on $H_0(\Sigma_g)$, $H_2(\Sigma_g)$, and $2g-2$ of the basis of $H_1(\Sigma_g)$; In the subspace for the remaining two basis of $H_1(\Sigma_g)$, their actions are represented by $\begin{pmatrix}
 1 & 1\\
& 1                                                                                                                                                                                                                                 
                                                                                                                                                                                     \end{pmatrix}$ and
$\begin{pmatrix}
1& \\
-1 &1  
 \end{pmatrix}
$ respectively. We then find the trace of the matrices $\begin{pmatrix}
                                             1&m\\ &1
                                            \end{pmatrix}
\cdot
\begin{pmatrix}
1 & \\
-n &1
\end{pmatrix}
$ and carry out the appropriate alternating sum.  In the end, the Lefschetz number is found to be $(2-2g+mn)$ - identical to the number of pairs of essential generators in this case.


The remaining arguments are similar.  The rank of the homology is at least the Eular characteristic, being the same as the Lefschetz number, but simultaneously cannot exceed the number of pairs of essential generators.  That leaves no other possibility than $2g-2+|mn|$.  So we have: 
\begin{theorem}

$HF^+(M(t_{\gamma}^mt_{\delta}^n), S_{g-2};\omega)=\A^{2g-2+|mn|},$ \;\;
$m\cdot n<0.$
\end{theorem}

Again, compare the result with Theorem 5.7 of Jabuka and Mark \cite{JM}; The ranks agree this time. 

\section{Multiple Dehn twists along a transverse pair of curves, case $m\cdot n>0$}

In this section, we compute the Heegaard Floer homology for the manifolds $M(t_{\gamma}^mt_{\delta}^n)$ where $m\cdot n>0$.  By symmetry, it is enough to consider the case $m, n>0$.   

\begin{figure}

\begin{center}
\includegraphics[width=3in]{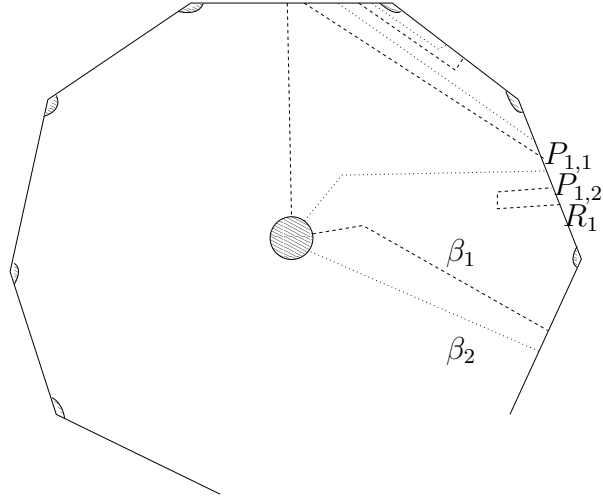}
\setlength{\unitlength}{1.02in}
\put(-0.7,1.2){$\beta_1$}
\put(-0.7,0.7){$\beta_2$}
\put(-0.2,1.7){$P_{1,1}$}
\put(-0.15,1.54){$P_{1,2}$}
\put(-0.1,1.38){$R_1$}
\caption{\label{HD3} The Heegaard Diagram for $M(t_\gamma^mt_\delta^n)$, when $m=1,n=2$.  $\beta_1$ is represented by the dashed curve, while $\beta_2$ is represented by the dotted curve.} 
\end{center}

\end{figure}

\begin{figure}

\begin{center}
\includegraphics[width=3in]{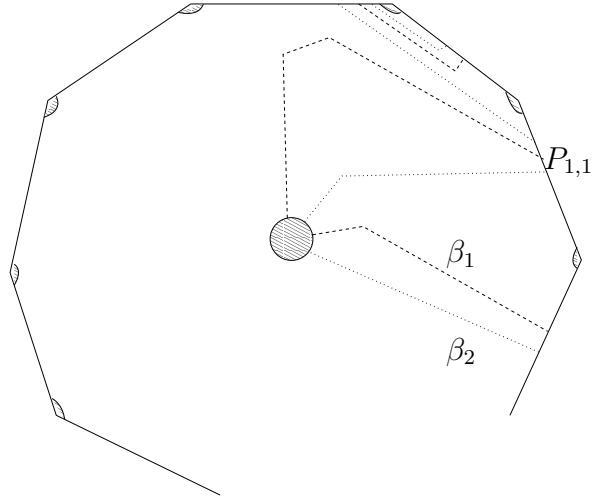}
\setlength{\unitlength}{1.02in}
\put(-0.7,1.2){$\beta_1$}
\put(-0.7,0.7){$\beta_2$}
\put(-0.2,1.68){$P_{1,1}$}
\caption{\label{HD4} The simplified Heegaard Diagram after isotopying on $\beta_1$. } 
\end{center}

\end{figure}

We have the following Heegaard diagram (Figure \ref{HD3}), and it can be subsequently simplified to Figure \ref{HD4} by an isotopy on $\beta_1$.  Note that the intersections $R_1$ and $P_{m,n}$ disappear there.  Within such a Heegaard diagram, there are $2(2g-4+mn)$ generators of $S_{g-2}$:
$$(A_{2g}, B_{2g}, L_1, R_2, \cdots,L_{2g-1})$$
$$...$$
$$(A_{2g}, B_{2g}, L_1,L_2, \cdots, R_{2g-1})$$
and \\
$$(A_{2g-1},B_{2g-1},L_1, R_2,\cdots, L_{2g})$$
$$...$$
$$(A_{2g-1},B_{2g-1},L_1,L_2 \cdots,R_{2g})$$
and \\
$$(A_{2g}, B_{2g}, P_{i,j}, L_2, \cdots, L_{2g-1}), (i,j)\neq (m,n)$$ 
$$(A_{2g-1},B_{2g-1}, P_{i,j},L_2, \cdots, L_{2g}), (i,j)\neq (m,n).$$

In the present case, it is possible and necessary to further partition all generators of $S_{g-2}$ according to their Chern classes.  Recall the first Chern class formula \cite[section 7.1]{OSzAnn2}:
$$\langle c_1(\s_y),[\mathcal{P}] \rangle =\chi(\mathcal{P})-2\overline{n}_z(\mathcal{P})+2\sum_{p\in y}\overline{n}_p(\mathcal{P}).$$
where $\s_y$ is a spin$^{c}$-structure corresponding to $y$.  Applying this formula, we find \\
$(A_{2g}, B_{2g}, P_{i,j}, L_2, \cdots, L_{2g-1}), (A_{2g-1},B_{2g-1}, P_{i,j},L_2, \cdots, L_{2g}), (i,j)\neq (m,n)$ lie on $mn-1$ different spin$^c$-structures, denoted by $\s_{i,j}$ respectively;  While all the remaining generators $(A_{2g}, B_{2g}, L_1, R_2, \cdots,L_{2g-1})$,...,$(A_{2g}, B_{2g}, L_1,L_2, \cdots, R_{2g-1})$, $(A_{2g-1},B_{2g-1},L_1, R_2,\cdots, L_{2g})$,...,
$(A_{2g-1},B_{2g-1},L_1,L_2 \cdots,R_{2g})$ lie on the other spin$^c$-structure, denoted by $\s_{m,n}$.

For each spin$^c$-structure $\s_{i,j}$, $(i,j)\neq (m,n)$, there are exactly two generators \\
$(A_{2g}, B_{2g}, P_{i,j}, L_2, \cdots, L_{2g-1})$, $(A_{2g-1},B_{2g-1}, P_{i,j}, L_2, \cdots, L_{2g})$ and an obvious holomorphic disk $D$ connecting them.  The argument from \cite[section 3]{ztwu} for 3-torus can thus be applied and shows  $$HF^+(M(t_{\gamma}^mt_{\delta}^n), \s_{i,j};\omega)=\A.$$

To determine the homology in the spin$^c$-structure $\s_{m,n}$, note its Euler characteristic is: 
\begin{align*}
 \chi(HF^+_{\s_{m,n}}) &=\tau_t(2g-2)-\sum_{(i,j)\neq (m,n)} \chi(HF^+_{\s_{i,j}}) \\
&=2-2g+mn-(mn-1) \\
&=3-2g.
\end{align*}

There are also exactly $2g-3$ pairs of essential generators, so we can repeat the argument within $\s_{m,n}$, as before, and conclude
$$HF^+(M(t_{\gamma}^mt_{\delta}^n), \s_{m,n};\omega)=\A^{2g-3}.$$

Putting everything together, we have:
\begin{theorem}

 $HF^+(M(t_{\gamma}^mt_{\delta}^n), S_{g-2};\omega)=\A^{2g-4+mn}, $ \;\;
$m\cdot n > 0.$

\end{theorem}

Again, compare the result with Theorem 5.3 of Jabuka and Mark \cite{JM}; our rank is smaller by two.  


\section{Multiple Dehn twists along a transverse pair of curves, case III}

\begin{figure}

\begin{center}
\includegraphics[width=4in]{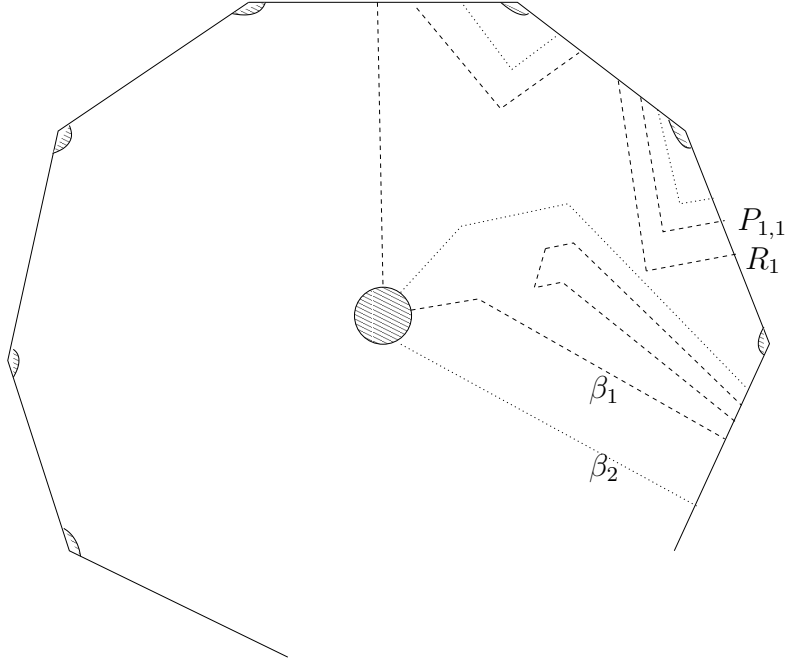}
\setlength{\unitlength}{1.36in}
\put(-0.7,1){$\beta_1$}
\put(-0.7,0.7){$\beta_2$}
\put(-0.1,1.5){$R_1$}
\put(-0.13,1.65){$P_{1,1}$}
\caption{\label{HD5} The Heegaard Diagram of $M(t_\gamma t_\delta t_\gamma)$.  An isotopy on $\beta_1$ can be carried out to cancel the pairs of intersection points $R_1$ and $P_{1,1}$. } 
\end{center}

\end{figure}

The manifolds considered here have the form $M(t_\gamma^{m_1}t_\delta^{n_1}t_\gamma^{m_2})$, where $m_1,m_2,n_1>0$.  The Heegaard diagram for the case $m_1=n_1=m_2=1$ is drawn in Figure \ref{HD5}.  An isotopy can be carried out for $\beta_1$ to remove the intersections $R_1$ and $P_{1,1}$.  Such an isotopy is available in general, so that we would have $2g-4+(m_1+m_2)n_1$ pairs of essential generators in the simplified Heegaard diagram.

The number of spin$^c$-structures in $S_{g-2}$ is $(m_1+m_2)n_1$, and these spin$^c$ structures are denoted by $\s_{i,j}$.  Similar to the case of section 5, there exists exactly one pair of essential generators in each $\s_{i,j}$ for $(i,j)\neq (m_1+m_2,n_1)$, and $2g-3$ pairs of essential generators in the remaining distinguished spin$^c$-structure $\s_{m_1+m_2,n_1}$.  

Hence, for all $(i,j)\neq (m_1+m_2, n_1)$, $$HF^+(M(t_\gamma^{m_1}t_\delta^{n_1}t_\gamma^{m_2}),\s_{i,j};\omega)=\A.$$ 
Since the Lefschetz number is $2-2g+(m_1+m_2)n_1$, we have: 
\begin{align*} \chi(HF^+_{\s_{m_1+m_2,n_1}}) &=\tau_t(2g-2)-\sum_{(i,j)\neq (m_1+m_2,n_1)} \chi(HF^+_{\s_{i,j}}) \\
&=2-2g+(m_1+m_2)n_1-((m_1+m_2)n_1-1) \\
&=3-2g
\end{align*}
The usual argument applies once more and shows: $$HF^+(M(t_\gamma^{m_1}t_\delta^{n_1}t_\gamma^{m_2}),\s_{m_1+m_2,n_1};\omega)=\A^{2g-3}.$$
Putting all the spin$^c$-structures together, we conclude with:
\begin{theorem}
 
$HF^+(M(t_\gamma^{m_1}t_\delta^{n_1}t_\gamma^{m_2}),S_{g-2};\omega)=\A^{2g-4+(m_1+m_2)n_1}.$

\end{theorem}

\section{Multiple Dehn twists along a transverse pair of curves, case IV}

Lastly, we consider the manifolds of the form $M(t_{\gamma}^{m_1}t_{\delta}^{n_1}\cdots t_{\gamma}^{m_k}t_{\delta}^{n_k})$, where $m_i\cdot n_j<0$.  In other words, they are the mapping tori of Dehn twists along $\gamma$ and $\delta$ with alternating signs.  

We compute the Lefschetz number.  If we denote the trace of the following matrix
$
\begin{pmatrix}
   1& m_1\\ 
&1
  \end{pmatrix}
\cdot
\begin{pmatrix}
 1 &\\ -n_1& 1
\end{pmatrix}
\cdots
\begin{pmatrix}
   1& m_k\\ 
&1
  \end{pmatrix}
\cdot
\begin{pmatrix}
 1 & \\ -n_k& 1
\end{pmatrix}
$ by $T$,
then the Lefschetz number is $4-2g-T$.  


Meanwhile, there are $2g+T-4$ pairs of essential generators.  This can be evidently seen by relating the intersections of $\alpha_i$ and $\beta_i$ to the trace of the matrix 
$
\begin{pmatrix}
   1& |m_1|\\ 
&1
  \end{pmatrix}
\cdot
\begin{pmatrix}
 1 &\\ |n_1|& 1
\end{pmatrix}
\cdots
\begin{pmatrix}
   1& |m_k|\\ 
&1
  \end{pmatrix}
\cdot
\begin{pmatrix}
 1 & \\ |n_k|& 1
\end{pmatrix}.
$

In our case, the number of pairs of essential generators is, once again, exactly the Lefschetz number.  Hence, we repeat our argument and conclude:
\begin{theorem}

$HF^+(M(t_{\gamma}^{m_1}t_{\delta}^{n_1}\cdots t_{\gamma}^{m_k}t_{\delta}^{n_k}), S_{g-2}; \omega)=\A^{2g-4+T}
$ \;\; $m_i\cdot n_j<0$.

\end{theorem}

\end{document}